\documentclass{amsart}

\newtheorem{thm}{Theorem}[section]
\newtheorem{lem}[thm]{Lemma}
\newtheorem{prop}[thm]{Proposition}

\theoremstyle{definition}

\theoremstyle{remark}

\numberwithin{equation}{section}

\newcommand{\C}{\mathbb{F}}
\newcommand{\Cc}{\mathcal{C}}
\newcommand{\Cs}{\Cc_\sigma}
\newcommand{\homcs}{\hom_{\Cc_\sigma}}
\newcommand{\tot}{\tilde{\otimes}}
\newcommand{\End}{\operatorname{End}}
\newcommand{\Endcs}{\End_{\Cs}}
\newcommand{\qhat}{\hat{q}}
\newcommand{\g}{\mathfrak{g}}
\newcommand{\h}{\mathfrak{h}}

\newcommand{\cqsln}{\mathbb{F}_q[SL(n)]}

\begin{document}

\author{Tim Hodges}
\title[Cremmer-Gervais $R$-matrices]{Generating functions for the coefficients of the Cremmer-Gervais $R$-matrices}
\address{University of Cincinnati, Cincinnati, OH 45221-0025,
U.S.A.}
\email{timothy.hodges@uc.edu}
\thanks{Supported in part by NSF grant DMS-9501484 and NSA grant MDA904-99-1-0026}
\date{March 18, 1999}

\begin{abstract}
	The coefficients of certain operators on $V\otimes V$ can be constructed  using generating functions. Necessary and sufficient conditions are given for some such operators to satisfy the Yang-Baxter equation. As a corollary we obtain a simple, direct proof that the Cremmer-Gervais $R$-matrices satisfy the Yang-Baxter equation. This approach also clarifies Cremmer and Gervais's original proof via the dynamical Yang-Baxter equation.
\end{abstract}
\maketitle

\section{Introduction}

	In \cite{BG}, Bilal and Gervais found a family of solutions to the dynamical Yang-Baxter equation (DYBE), or Gervais-Neveu equation. By a change of basis argument they found an extremely interesting family of solutions to the usual Yang-Baxter equation, now known as the Cremmer-Gervais $R$-matrices. This work raised questions about the connections between the DYBE and the YBE which have yet to be addressed. It is also natural to ask whether there are other, more direct approaches which lead to a proof that the Cremmer-Gervais $R$-matrices satisfy the YBE. A direct proof was given by the author in \cite{H} but it involved some rather lengthy technical calculations. The work of Etingof and Kazhdan \cite{EK} yields a quantization of solutions of the classical Yang-Baxter equation which in theory yields the Cremmer-Gervais $R$-matrices as a special case. However it has not yet been possible to make this construction explicit. 

	We show here that the coefficicients of the Cremmer-Gervais $R$-matrices can be realised as the coefficients of certain simple ``generating functions''. For operators of such form we give necessary and sufficient conditions for the YBE to hold. These conditions are then easy to verify for the Cremmer-Gervais $R$-matrices. As a further example of operators of such form, we consider the operator $\eta$ used in \cite{H}. Again it follows easily from the theorem that $\eta$ satisfies the YBE.

	Next we turn to discussing the original argument in \cite{CG}. We begin by recalling briefly the axiomatic formulation of the DYBE given by Etingof and Varchenko \cite{EV} and the standard solution found in \cite{BG}. We note at this stage that the coefficients of these matrices are expressible using the same functions that arose earlier as components of the generating functions. We then prove a general change of basis result which enables one to pass from solutions of the DYBE to solutions of the YBE. Finaly we show how the appropriate change of basis yields the Cremmer-Gervais $R$-matrices. Again the use of generating functions clarifies and simplifies these calculations.

	The author would like to thank Jintai Ding for his comments and suggestions.

\section{Generating functions and the Yang-Baxter equation}

	Let $V$ be an $n$-dimensional vector space over a field $\C$ with basis $\{e_1, e_2, \dots, e_n\}$. An operator $\gamma \in \End(V\otimes V)$ given by
$$\gamma(e_i \otimes e_j) = \sum_{k,l}\gamma^{kl}_{ij} e_k \otimes e_l$$
will be said to be homogeneous if $\gamma^{kl}_{ij}\neq 0$ only if $i+j=k+l$. In this case we may write $\gamma(i,j,k) = \gamma^{k, i+j-k}_{ij}$ so that
$$\gamma(e_i \otimes e_j) = \sum_k \gamma(i,j,k) e_k \otimes e_{i+j-k}.$$
We will call the polynomials 
$G^\gamma_{i,j}(x) =\sum_k \gamma(i,j,k)x^k$ the {\em generating functions} for the coefficients of $\gamma$.

	As a simple example, take $\eta$ to be the homogeneous operator defined by
$$
\eta(i,j,k) = \begin{cases}
	1 & \text{ if $i\leq k < j$} \\
	-1 & \text{ if $j\leq k < i$} \\
0 & \text{ otherwise }
	\end{cases}  
$$
Then it is easy to see that the generating functions for $\eta$ are
$$
	\sum_k \eta (i,j,k)x^k = \frac{x^i - x^j}{1-x}
$$

	The Cremmer-Gervais operators in their most general two-parameter form are
$$
	\rho_p (e_i \otimes e_j)= q p^{i-j} e_j \otimes e_i + \sum_k \qhat p^{i-k} \eta(i,j,k) e_k \otimes e_{i+j-k}
$$
where $q$ and $p$ are non-zero elements of the base field $\C$ and $\qhat = q -q^{-1}$ \cite{CG,H2}.

\begin{lem}\label{rho1}
	The generating functions for the Cremmer Gervais operator $\rho_p$ are
$$
G^{\rho_p}_{i,j}(x) = \frac{\qhat}{1-p^{-1}x}x^i +
 p^{i-j}\frac{q^{-1}-qp^{-1}x}{1-p^{-1}x}x^j.
$$
In particular in the case $p=1$, they become
$$
G^{\rho_1}_{i,j}(x)  = \frac{\qhat}{1-x}x^i + \frac{q^{-1}-qx}{1-x}x^j.
$$
\end{lem}

\begin{proof}
\begin{align*}
\sum_k \rho_p&(i,j,k) x^k  = q p^{i-j} x^j + \sum_k \qhat p^{i-k} \eta(i,j,k)x^k\\
&= q p^{i-j} x^j + \sum_k \qhat p^{i} \eta(i,j,k)(p^{-1}x)^k\\
&= q p^{i-j} x^j + \qhat p^i \frac{(p^{-1}x)^i -(p^{-1}x)^j}{1 - p^{-1}x}\\
&= \frac{\qhat}{1-p^{-1}x}x^i +
 p^{i-j}\frac{q^{-1}-qp^{-1}x}{1-p^{-1}x}x^j.
\end{align*}
\end{proof}

	Thus the generating functions for both $\eta$ and $\rho_1$ are of the form
$\alpha(x)x^i + \beta(x)x^j$ where $\alpha(x)$ and $\beta(x)$ are rational functions independent of $i$ and $j$. We now determine general conditions on $\alpha(x)$ and $\beta(x)$ which guarantee that such an operator satisfies the Yang-Baxter equation. Here we shall be using the ``braid version'' of the Yang-Baxter equation: $R_{12}R_{23}R_{12} = R_{23}R_{12}R_{23}$. Note that we do not require that a solution of the Yang-Baxter be invertible.

\begin{thm}\label{gfybe}
	Let $\gamma  \in \End(V\otimes V)$ be a homogeneous operator for which the generating functions are of the form $\alpha(x)x^i + \beta(x)x^j$. Then $\gamma$ satisfies the Yang-Baxter equation if and only if $\beta(x)=0$ or the following two identities hold:
\begin{enumerate}
\item $\alpha(x)\alpha(y) = \alpha(xy^{-1})\alpha(y) + \alpha(x)\alpha(yx^{-1})$
\item $\alpha(xy^{-1}) ^2\alpha(y) + \beta(xy^{-1}) \alpha(x) \beta(yx^{-1})
	= \alpha(y)^2 \alpha(xy^{-1}) + \beta(y)\alpha(x)\beta(y^{-1})$
\end{enumerate}
\end{thm}

\begin{proof}
	If $\beta(x)=0$ then $\alpha(x)$ must be a scalar and the YBE trivially holds. So assume that $\beta(x)\neq 0$.

	By applying both sides of the equation $\gamma_{23}\gamma_{12}\gamma_{23}=\gamma_{12}\gamma_{23}\gamma_{12}$ to $e_i \otimes e_j \otimes e_k$ and comparing coefficients one sees that a homogeneous operator $\gamma$ will satisfy the Yang Baxter equation if and only if 
\begin{multline*}
\sum_a \gamma(j,k,a) \gamma(i,a,c)\gamma(i+a-c,j+k-a,h)\\
= \sum_s \gamma(i,j,s)\gamma(i+j-s,k,h+c-s) \gamma(s,h+c-s,c)
\end{multline*}
for all $i$, $j$, $k$, $h$ and $c$. 
Equivalently, $\gamma$ will satisfy the Yang-Baxter equation if
\begin{multline*}
\sum_{a,c,h} \gamma(j,k,a) \gamma(i,a,c)\gamma(i+a-c,j+k-a,h)x^c y^h\\
= \sum_{s,c,h} \gamma(i,j,s)\gamma(i+j-s,k,h+c-s) \gamma(s,h+c-s,c)x^c y^h
\end{multline*}
for all $i$, $j$ and $k$. Now 
\begin{align*}
&\sum_{a,c,h} \gamma(j,k,a) \gamma(i,a,c)\gamma(i+a-c,j+k-a,h)x^c y^h\\
&= \sum_{a,c} \gamma(j,k,a) \gamma(i,a,c)(\alpha(y) y^{i+a-c} + \beta(y)y^{j+k-a})x^c\\
&= \sum_{a,c} \gamma(j,k,a) \gamma(i,a,c)[\alpha(y) y^{i+a}(xy^{-1})^c + \beta(y)y^{j+k-a}x^c]\\
&= \sum_{a} \gamma(j,k,a) 
	[\alpha(y) y^{i+a}[\alpha(xy^{-1})(xy^{-1})^i + \beta(xy^{-1}) (xy^{-1})^a]\\
 & \qquad + \beta(y)y^{j+k-a}[\alpha(x) x^i + \beta(x)x^a]]
\end{align*}
\begin{align*}
&= \sum_{a} \gamma(j,k,a) 
	[\alpha(y)\alpha(xy^{-1})x^i y^a + \alpha(y)\beta(xy^{-1}) y^i x^a\\
 & \qquad + \beta(y)\alpha(x) y^{j+k} x^i y^{-a} + \beta(y)\beta(x)y^{j+k}(xy^{-1})^a]\\
&=x^iy^k[\beta(y)\alpha(x)\alpha(y^{-1}) + \alpha(y)\alpha(xy^{-1})\beta(y)]\\
&\qquad+x^iy^j[\alpha(y)\alpha(xy^{-1}) \alpha(y) + \beta(y)\alpha(x) \beta(y^{-1})]\\
&\qquad + x^jy^i \alpha(y) \beta(xy^{-1})\alpha(x) + x^ky^i\alpha(y)\beta(xy^{-1})\beta(x)\\
&\qquad + x^jy^k\beta(y) \beta(x) \alpha(xy^{-1}) + x^ky^j \beta(y) \beta(x)\beta(xy^{-1})
\end{align*}
Similarly,
\begin{align*}
 & \sum_{s,c,h} \gamma(i,j,s)\gamma(i+j-s,k,h+c-s) \gamma(s,h+c-s,c)x^c y^h\\
& = x^iy^j[\alpha(xy^{-1})^2\alpha(y) + \beta(xy^{-1})\alpha(x)\beta(yx^{-1})]\\
& \qquad + x^j y^i[\alpha(xy^{-1})\alpha(y)\beta(xy^{-1}) + \beta(xy^{-1})\alpha(x)\alpha(yx^{-1})]\\
&\qquad + x^iy^k \alpha(xy^{-1})\beta(y)\alpha(x) + x^ky^i\alpha(y)\beta(xy^{-1})\beta(x)\\
&\qquad + x^jy^k\beta(y) \beta(x) \alpha(xy^{-1}) + x^ky^j \beta(y) \beta(x)\beta(xy^{-1})
\end{align*}
Comparing coefficients then yields the theorem.
\end{proof}

	Since very few operators do have generating functions of the simple form $\alpha(x)x^i + \beta(x)x^j$ one can completely classify the solutions of the Yang-Baxter equation that arise in this form. They are precisely the solutions found in \cite{H} together with their ``transposes''.
	Notice that the ``flip'' operator $P(e_i \otimes e_j) = e_j \otimes e_i$ has generating functions $G^P_{i,j}(x) = x^j$ and for the identity $I$, $G^I_{i,j}(x) = x^i$. 

\begin{prop} The homogeneous operators having generating functions of the form $\alpha(x)x^i + \beta(x)x^j$ for $\alpha(x), \beta(x) \in \C(x)$ are the operators of the form $\gamma= aI + bP + c \eta$. For these operators the generating functions are:
$$
	G^\gamma_{i,j}(x) = \left( a + \frac{c}{1-x}\right) x^i + \left( b - \frac{c}{1-x}\right) x^j.
$$
\end{prop}

\begin{proof} First note that because $I$, $P$ and $\eta$ have generating functions of the desired form, so does any linear combination. Hence it remains to show that these are the only possibilities.

	Let $\alpha(x), \beta(x) \in \C(x)$ and suppose that $\alpha(x)x^i + \beta(x)x^j$ is a polynomial of degree less than or equal to $n$ for all $1 \leq i,j\leq n$. It is easily checked that if $\alpha(x), \beta(x) \in \C[x,x^{-1}]$ then they must be scalar. Embed the rational function field $\C(x)$ into the Laurent power series ring $\C[[x,x^{-1}]]$ in the usual way. Considering the cases $i=j=n$ and $i=n$, $j=n-1$ yields that $\alpha(x)$ and $\beta(x)$ have eventually constant coefficients of opposite sign. Hence,
$$
\alpha(x) = \alpha'(x) + \frac{c}{1-x}; \quad \beta(x) = \beta'(x) - \frac{c}{1-x}
$$
where $\alpha'(x),\beta'(x) \in \C[x,x^{-1}]$ and $c\in \C$. But then the functions 
$$
\alpha(x)x^i + \beta(x)x^j - P^{c\eta}_{i,j}(x) =\alpha'(x)x^i + \beta'(x)x^j
$$ are also all polynomial. Hence $\alpha'(x)$ and $\beta'(x)$ must be scalars, as required.
\end{proof}

\begin{thm}
	The homogeneous operators $\gamma$ having generating functions of the form $\alpha(x)x^i + \beta(x)x^j$ for $\alpha(x), \beta(x) \in \C(x)$ and which satisfy the Yang-Baxter equation are of one of the two following forms.
\begin{enumerate}
\item $aP + b \eta$, for $a, b \in \C$.
\item $aP + b(I-\eta)$ for $a,b \in \C$.
\end{enumerate}
\end{thm}

\begin{proof} We need to determine which pairs of functions of the form
$$
\alpha(x) = a + \frac{c}{1-x}; \quad \beta(x) = b - \frac{c}{1-x}
$$
satisfy the equations
\begin{enumerate}
\item $\alpha(x)\alpha(y) = \alpha(xy^{-1})\alpha(y) + \alpha(x)\alpha(yx^{-1})$
\item $\alpha(xy^{-1}) ^2\alpha(y) + \beta(xy^{-1}) \alpha(x) \beta(yx^{-1})
	= \alpha(y)^2 \alpha(xy^{-1}) + \beta(y)\alpha(x)\beta(y^{-1})$
\end{enumerate}

	The first equation is satisfied if 
$$
\alpha(x) = \frac{c}{1-x} \quad \text{or} \quad \alpha(x) = \frac{cx}{1-x}.
$$
On the other hand $\alpha(x)$ cannot have any non-zero roots. For if $\alpha(d) =0$, then $\alpha(dy^{-1})\alpha(y) = 0$ which is impossible if $d\neq 0$. So the above are indeed the only possibilities for $\alpha(x)$.

	Now suppose that $\alpha(x)=c/(1-x)$ and let $\beta(x) = b - \alpha(x)$. Then
\begin{align*}
\beta(x)\beta(x^{-1}) &= (b- \alpha(x))(b-\alpha(x^{-1}))\\
		&= b^2 - b(\alpha(x) + \alpha(x^{-1})) +\alpha(x)\alpha(x^{-1})\\
		&= b(b-1) + \alpha(x)\alpha(x^{-1})
\end{align*}
Using this identity, the second equation follows easily from the first equation. This yields all solutions of the first form. A similar analysis of the case $\alpha(x)=cx/(1-x)$ yields all solutions of the second form. The restriction on $\alpha(x)$ implies that these are the only possibilities.
\end{proof}

\noindent {\em Remarks}

\begin{enumerate}
\item The Cremmer-Gervais operator $\rho_1$ is the special case $qP + \qhat \eta$. The fact that the more general Cremmer-Gervais operators satisfy the Yang-Baxter equation can be deduced using some elementary twisting arguments.
\item The fact that the operators of the first type satisfy the Yang-Baxter equation was proved in \cite{H} using some rather complex and unilluminating identities for the functions $\eta(i,j,k)$. Almost all  the identities proved there can be explained with the use of generating functions.
\item It is well-known that if an operator $\gamma$ satisfies the Yang-Baxter equation then so does its ``transpose'' $P\gamma P$. The operators of the second type above are precisely the transposes of those of the first type.
\item Form the identities $\eta^2 = \eta$, $\eta P = - \eta$ and $P \eta = \eta + P +I$, it follows easily that $\gamma= aP + b \eta$ satisfies
$$(\gamma-a)(\gamma+(a-b)) =0.
$$
Hence $\gamma$ is invertible if and only if $a(a-b)\neq 0$ and satisfies the Hecke relation $(\gamma - q)(\gamma+q^{-1})=0$ when $a=q$ and $b = \qhat$ \cite{H}.
\item	The polynomial functions 
$$\sum_{a,c,h} \gamma(j,k,a) \gamma(i,a,c)\gamma(i+a-c,j+k-a,h)x^c y^h$$
are generating functions for the homogeneous operator $\gamma_{23}\gamma_{12}\gamma_{23} \in \End (V \otimes V\otimes V)$. In the case when $\gamma=\eta$ they have the particularly simple form
$$
\frac{x^ky^j + x^j y^i + x^i y^k - x^k y^i - x^i y^j -x^j y^k}{(x-1)(y-1)(xy^{-1}-1)}.
$$
\end{enumerate}

\section{The DYBE and the proof of Cremmer and Gervais}

	We now briefly present the original proof of Cremmer and Gervais in an axiomatic algebraic framework. We begin by discussing the dynamical Yang-Baxter equation. Our approach is essentially that of Etingof and Varchenko \cite{EV,EV2}.

\subsection{The tensor category $\Cs$ and the $\sigma$-DYBE}

	Let $H$ be a commutative cocommutative Hopf algebra. Let $B$ be an $H$-module algebra with structure map,
$$
	\sigma: H \otimes B \to B
$$
Denote by $\Cc$ the category of right $H$-comodules. Define a new category $\Cs$ whose objects are right $H$-comodules but whose morphisms are
$$
\homcs (V, W) = \hom_H (V, W \otimes B)
$$ where $B$ is given a trivial comodule structure. Composition of morphisms is given by the natural embedding of $\hom_H (V, W \otimes B)$ inside $\hom_H (V\otimes B, W \otimes B)$.

	We now define a tensor product on this category. Define a bifunctor
$$
\tot : \Cs \times \Cs \to \Cs.
$$
For objects $V$ and $W$, $V\tot W$ is the usual tensor product of $H$ comodules $V \otimes W$. In order to define the tensor product of two morphisms, notice first that we can define, for any $H$-comodule $W$,  a linear twist map $\tau : B \otimes W \to W \otimes B$ by
$$
	\tau(b \otimes w) = w_0 \otimes \sigma(w_1 \otimes b).
$$
Then for any pair of morphisms $f : V \to V'$ and $g : W \to W'$, define
$$
f \tot g = (1 \otimes m_B)(1 \otimes \tau \otimes 1)(f \otimes g)
$$

\begin{thm}\cite{EV,EV2} The bifunctor $\tot$  makes $\Cs$ into a tensor category.
\end{thm}

	Let $V \in \Cs$ For any $R \in \Endcs(V \tot V)$ we define elements of $\Endcs (V \tot V \tot V)$, $R_{12} = R \tot 1$ and $R_{23} = 1 \tot R$. Then $R$ is said to satisfy the $\sigma$-dynamical Yang-Baxter equation ($\sigma$-DYBE) if
$$
 	R_{12}R_{23}R_{12} = R_{23}R_{12}R_{23}.
$$

	A more traditional formulation of the dynamical Yang-Baxter equation is the following. Let $\h$ be a Cartan subalgebra of a semisimple Lie algebra $\g$ and let $V$ be a module over $\g$.  Consider a meromorphic function
$$
	R: \h^* \to \End (V \otimes V)_, \quad \lambda \mapsto R(\lambda)
$$
Define $R_{23}(\lambda)$ in the usual way but define $R_{12}(\lambda + h^{(3)})$ by
$$
	R_{12}(\lambda + h^{(3)})(u\otimes v \otimes w) = R(\lambda + \mu)(u \otimes v)\otimes w
$$
if $w$ is a weight vector of weight $\mu$. The dynamical Yang-Baxter equation (in its braided form) is then
$$
	R_{12}(\lambda+h^{(3)}) R_{23}(\lambda)R_{12}(\lambda+h^{(3)})= R_{23}(\lambda)R_{12}(\lambda+h^{(3)}) R_{23}(\lambda)
$$
	
\subsection{The standard solution}\label{stansol}

	We now describe the standard (Bilal-Gervais) solution to this equation. 

	Let $\cqsln$ be the usual quantum group and let $T$ be the usual maximal torus of $SL(n)$. Then $H=\C[T]$ is a homomorphic image of $\cqsln$. Alternatively we may think of $H$ as the group algebra of the weight lattice $P$, so $H = \C[K_\lambda \mid \lambda \in P]$. Let $V$ be the standard comodule over $\cqsln$. Then $V$ has a basis $\{e_i\}$ of weight vectors with weights $\nu_i$. Denote the structure map by $\rho : V \to V \otimes \C[T]$. Then $\rho(e_i) = e_i \otimes K_{\nu_i}$.
	
	Now define $B$ to be the field of fractions of the subalgebra of $H$ generated by the root lattice $Q \subset P$. That is,
$$
	B = \text{Frac} (\C[ K_\alpha \mid \alpha \in Q])
$$
Define an action of $\sigma: H \otimes B \to B$ by
$$
	\sigma(K_\lambda \otimes K_\alpha) = q^{(\lambda, \alpha)} K_\alpha.
$$ 
Denote $\sigma (K_\lambda\otimes b)$ by $b^{\lambda}$. Let $\alpha_{ij}$ be the usual root $\nu_i-\nu_j$. Henceforth we need to assume that $q^2\neq 1$.

\begin{thm} \cite{BG,CG} \label{ssdybe}
	Suppose that $R \in \Endcs (V \tot V)$ is given by
$$
	R(e_i \otimes e_j) =  e_i \otimes e_j \otimes  \alpha((K_{\alpha_{ji}}q^{-\delta_{ij}})^2)
	+e_j \otimes e_i \otimes \beta((K_{\alpha_{ji}}q^{-\delta_{ij}})^2) 
$$
where $\alpha(x)=(q-q^{-1})/(1-x)$ and $\beta(x) =(q^{-1}-qx)/(1-x)$.
Then $R$ satisfies the $\sigma$-DYBE
\end{thm}

	This operator is the reformulation from \cite{CG} of the solution given  in \cite{BG}. Many other proofs now exist, both direct and conceptual \cite{ABRR,EV,EV2,I,JKOS}. Usually the base field is assumed to be the complex numbers but the more direct proofs work over any field.

\subsection{Change of basis and the YBE}

	Suppose that $f : V \to V'$ and $g: W \to W'$ are {\em linear } maps between $H$-comodules. Then we may again form $f \tot g$ as above. In this situation we no longer have that 
$$
	(f \tot g)(f' \tot g') = ff' \tot gg'
$$
because in general $(1 \tot g)(f \tot 1) \neq f \tot g$. However the following identities remain true:
\begin{enumerate}
\item $(f \tot 1)(f'\tot 1) = ff'\tot 1$
\item $(1 \tot g)(1 \tot g') = 1\tot gg'$
\item $(f \tot 1) (1 \tot g) = f \tot g$.
\end{enumerate}
Moreover $(1 \tot g)(f \tot 1) = f \tot g$ if $g$ is a comodule homomophism or if $f(V) \subset V'\otimes \C$.

 	Now let $A : V \to V \otimes B$ be a linear map. Define $A_1 = A \tot 1$, $A_2 =1 \tot A$ from $V \otimes V$ to $V \otimes V \otimes B$ and similarly $A_1$, $A_2$ and $A_3$ from $V \otimes V \otimes V$ to $V \otimes V \otimes V \otimes B$.

\begin{prop}\label{cob}
	Let $R \in \Endcs(V \tot V)$ be a solution of the $\sigma$-DYBE. Let $A: V \to V \otimes B$ be a linear map. Set $R^A = A_2^{-1} A_1^{-1} R A_1 A_2$. Suppose that
$ R^A(V \otimes V) \subset  V \otimes V\otimes \C$.
Then $R^A$ satisfies the Yang-Baxter equation.
\end{prop}

\begin{proof}
Clearly,
$$
	R^A_{23}R^A_{12}R^A_{23} = A_3A_2^{-1}R_{23}A_2 A_3 A_2^{-1} A_1^{-1} R_{12} A_1 A_2 
	A_3^{-1} A_2^{-1} R_{23}A_2 A_3
$$
and
$$
	R^A_{12} R^A_{23} R^A_{12} = A_2^{-1} A_1^{-1} R_{12} A_1 A_2 A_3^{-1}A_2^{-1} R_{23}
	A_2 A_3 A_2^{-1} A_1^{-1} R_{12} A_1 A_2
$$
Now $A_1$ commutes with $R_{23}$ because $R_{23}$ is a comodule morphism and $A_3$ commutes with  $R^A_{12}$ by the hypothesis. Using these facts we can see that
$$
	R^A_{12}R^A_{23}R^A_{12}= A_3^{-1} A_2^{-1} A_1^{-1} R_{12} R_{23} R_{12}A_1 A_2 A_3
$$
and
$$
	R^A_{23}R^A_{12}R^A_{23}= A_3^{-1} A_2^{-1} A_1^{-1} R_{23} R_{12} R_{23}A_1 A_2 A_3
$$
	Hence
$$
	R^A_{12}R^A_{23}R^A_{12}=R^A_{23}R^A_{12}R^A_{23}
$$
	Since $R^A$ is a map from $V \otimes V $ to $V \otimes V$, it therefore satisfies the usual Yang-Baxter equation.

\end{proof}

\subsection{The Cremmer-Gervais $R$-matrix}

We now return to the set-up of section \ref{stansol}. In particular, $R$ will denote the standard solution of the $\sigma$-DYBE given in Theorem \ref{ssdybe},
$$
	R(e_i \otimes e_j) =  e_i \otimes e_j \otimes  \alpha((K_{\alpha_{ji}}q^{-\delta_{ij}})^2)
	+e_j \otimes e_i \otimes \beta((K_{\alpha_{ji}}q^{-\delta_{ij}})^2) 
$$
The standard Cremmer-Gervais $R$-matrix is the operator $\rho=\rho_{q^{2/n}}$ described above.  That is
$$
\rho(e_i \otimes e _j) = qq^{2(i-j)/n} e_j \otimes e_i + \sum_k \qhat q^{2(i-k)/n} \eta(i,j,k) e_k \otimes e_{i+j-k}
$$
We now deduce that $\rho$ is a solution of the Yang-Baxter equation by showing that it can be obtained from $R$ by an appropriate ``change of basis''. 

\begin{lem} \label{cgcob}
	Let $A: V \to V \otimes B$ be the linear map 
$A(e_i) = \sum_a e_a \otimes K_{\nu_a}^{-2i}$. Then $R^A = \rho$.
\end{lem}

\begin{proof}
	We prove that $R A_1 A_2 = A_1 A_2 \rho$. In matrix form this is equivalent to
$$
\sum_{c,d} R^{ms}_{cd} (A^c_i)^{\nu_d} A_j^d = \sum_{k,l} \rho^{kl}_{ij} (A^m_k)^{\nu_s} A_l^s.
$$
From Lemma \ref{rho1} we have that
$$
\sum_k \rho(i,j,k)q^{2k/n}x^k = q^{2i/n}[\alpha(x)x^i + \beta(x)x^j]
$$
Note also that (abreviating $K_{\nu_m}$ by $K_m$)
$$
A^m_a = K^{-2a}_m \quad \text{and} \quad (A_a^m)^{\nu_s} = q^{-2a\delta_{ms} +2a/n}A^m_a
$$
Hence 
\begin{align*}
\text{RHS} &= \sum_a \rho(i,j,k)(K_m^{-2a})^{\nu_s}K_s^{-2(i+j-a)}\\
	&= K_s^{-2(i+j)}\sum_a \rho(i,j,k)q^{2a/n}(K_sK_m^{-1}q^{-\delta_{ms}})^{2a}\\
	&= K_s^{-2(i+j)}q^{2i/n}[\alpha((K_sK_m^{-1}q^{-\delta_{ms}})^2)(K_sK_m^{-1}q^{-\delta_{ms}})^{2i} +\\
	&\qquad \quad \beta((K_sK_m^{-1}q^{-\delta_{ms}})^2)(K_sK_m^{-1}q^{-\delta_{ms}})^{2j}]\\
	&= \alpha((K_sK_m^{-1}q^{-\delta_{ms}})^2)(K_m^{-2i})^{\nu_s}K_s^{-2j}
	+ \beta((K_sK_m^{-1}q^{-\delta_{ms}})^2)(K_s^{-2i})^{\nu_m}K_m^{-2j}\\
	&= \sum_{c,d}R^{ms}_{cd} (A_i^c)^{\nu_d} A_j^d
\end{align*}
as required.
\end{proof}

\begin{thm} \cite{CG} The operator $\rho$ satisfies the Yang-Baxter equation.
\end{thm}

\begin{proof} The result follows from Lemma \ref{cgcob} and Proposition \ref{cob}.
\end{proof}

\end{document}